\documentclass[aip,jcp,onecolumn]{revtex4-1}
\usepackage{setspace}
\usepackage{epstopdf}
\usepackage{graphicx}
\usepackage{subfigure}
\usepackage{mathtools}
\usepackage{amsmath}
\usepackage{amssymb}

\usepackage{float}

\def\fatn{\mathbf{n}}
\def\fatr{\mathbf{r}}

\def\fatx{\mathbf{x}}

\def\faty{\mathbf{y}}
\def\fatY{\mathbf{Y}}

\def\fatp{\mathbf{p}}

\usepackage{newfloat,algcompatible}
\usepackage[size=small]{caption}
\usepackage{etoolbox}

\AtBeginEnvironment{algorithm}{\noindent\par\nobreak\vskip-5pt}
\usepackage{newfloat}

\DeclareFloatingEnvironment[
    fileext=loa,
    listname=List of Algorithms,
    name=ALGORITHM,
    placement=tbhp,
]{algorithm}
\DeclareCaptionFormat{algorithms}{\vskip-15pt\hrulefill\par#1#2#3\vskip-6pt\hrulefill}
\begin{document}
\singlespacing
\title{Single molecule simulations in complex geometries with embedded dynamic one-dimensional structures}
\author{Stefan Hellander}
\affiliation{Department of Information Technology, Uppsala University, Box 337, SE-75105, Uppsala, Sweden.} 
\email{stefan.hellander@it.uu.se}
\begin{abstract}
Stochastic models of reaction-diffusion systems are important for the study of biochemical reaction networks where species are present in low copy numbers or if reactions are highly diffusion limited. In living cells many such systems include reactions and transport on one-dimensional structures, such as DNA and microtubules. The cytoskeleton is a dynamic structure where individual fibers move, grow and shrink. In this paper we present a simulation algorithm that combines single molecule simulations in three-dimensional space with single molecule simulations on one-dimensional structures of arbitrary shape. Molecules diffuse and react with each other in space, they associate to and dissociate from one-dimensional structures as well as diffuse and react with each other on the one-dimensional structure. A general curve embedded in space can be approximated by a piecewise linear curve to arbitrary accuracy. The resulting algorithm is hence very flexible. Molecules bound to a curve can move by pure diffusion or via active transport, and the curve can move in space as well as grow and shrink. The flexibility and accuracy of the algorithm is demonstrated in four numerical examples.
\end{abstract}
\pacs{82.37.Np,82.39.-k,87.10.Mn,87.16.Ka}
\maketitle
\section{Introduction}

It has been established that modeling biochemical reaction networks with a macroscopic deterministic model, for instance with ordinary differential equations (ODEs), can be unsuitable when some species are present in low copy numbers or if reactions are highly diffusion limited \cite{MAA2,ElE}. Such systems may only be accurately modeled by a discrete stochastic model. Usually two levels of stochastic modeling are considered: the mesocopic level and the microscopic level.

The mesoscopic scale is described by the chemical master equation (CME) which is a discrete, stochastic model. Molecules are treated as discrete entities but assumed to be well stirred with a uniform distribution within the simulation volume. The copy number of each species is computed but the position of individual molecules is not tracked. Exact trajectories of such a system can be generated with the stochastic simulation algorithm (SSA) by Gillespie \cite{gillespie}.

Systems for which the condition of spatial homogeneity does not hold can be modeled at the mesoscopic level by subdividing the domain into smaller compartments, where in each compartment the well stirred condition is assumed to hold. Diffusion can now be modeled as discrete jumps between the compartments. The governing equation of the system is the reaction-diffusion master equation (RDME) and exact trajectories can be computed with the next subvolume method (NSM) \cite{ElE}, which is an efficient implementation of the SSA in a spatial setting. Several software packages have been developed for simulating systems at the mesoscopic level, including MesoRD \cite{mesoRD}, URDME \cite{DrEnHe} and STEPS \cite{WilsSchutter}. The spatial accuracy of the RDME depends on the size of the voxels. If the voxels are too large we might not resolve the spatial dynamics of the system accurately, but on the other hand it has also been shown that when the size of the voxels becomes too small the accuracy deteriorates \cite{Isaacson2,HeHePe,ErbChap}. By allowing the reaction rates in the system to depend on the size of the voxels in a general way, accuracy can be retained, but at most up to a critical lower bound of the voxel size, at which point the model breaks down regardless of how the reaction rates are chosen \cite{ErbChap,FBSE10,HeHePe}. The mesoscopic model may therefore be unsuitable for systems where a very high spatial accuracy is required to capture important dynamics. 

The microscopic level, or microscale, usually refers to models where the positions of single molecules are tracked. Either we compute the continous position in space or let the individual molecules move on a lattice, as in Ref. \onlinecite{TaArTo}. Reactions occur with some probability per unit time when molecules collide or when they are within some prescribed reaction radius. There are two models that have been used extensively to model systems at the microscale: the Smoluchowski model \cite{Smol} and the Doi model \cite{doi1,doi2}. In this paper we will only be concerned with the Smoluchowski model.

In Ref. \onlinecite{TaTNWo10} they show that fine-grained spatial correlations can affect the macroscopic behaviour of a system. Furthermore, the dynamics of the system cannot be captured on the macroscopic or mesoscopic level but rather require a microscopic level of simulation. This is due to highly diffusion-limited reactions and fast rebindings that the mesoscopic model cannot resolve without hitting the critical lower bound of the voxel size. This motivates development of efficient methods for simulating systems at the microscopic level.

The Green's function reaction dynamics (GFRD) algorithm \cite{ZoWo5a} is an efficient method to generate accurate trajectories of the Smoluchowski model. A full $N$-body problem is divided into many $1$- and $2$-body problems that can be simulated independently by choosing the time step in the algorithm sufficiently small. For these subproblems the Smoluchowski equation can be solved analytically or numerically. Improvements of GFRD have been developed in e.g. Refs. \onlinecite{SHeLo11,TaTNWo10}. Sampling random numbers from the analytical solution of the Smoluchowski equation is complicated and expensive, but in Ref. \onlinecite{SHeLo11} this step is simplified by solving the Smoluchowski equation using first order or second order operator splitting schemes. In Ref. \onlinecite{TaTNWo10} the algorithm is made exact in both space and time by combining GFRD with the first-passage kinetic Monte Carlo (FPKMC) algorithm developed in Ref. \onlinecite{FPKMC}, and this improved algorithm is called eGFRD. The main idea in eGFRD is to introduce protective spheres around single molecules and pairs of molecules and then sample exact times for the molecules to exit the spheres.

Software that implements eGFRD is eCell \cite{ecell}. Other software packages that simulate systems with the Smoluchowski model as the target model are Smoldyn \cite{AnB}, MCell \cite{MCell08} and Spatiocyte \cite{TaArTo}. Both Smoldyn and MCell discretize time while GFRD is exact in time, thus making them less accurate. On the other hand, Smoldyn and MCell can be significantly more efficient for larger systems. In Spatiocyte molecules move on a lattice of densely packed spheres.

Instead of simulating the system with a purely microscopic method one can combine the microscopic level and the mesoscopic level in a hybrid method \cite{HeHeLo,KlArHe,FlChEr11}. By restricting the microscale simulations to the species or the parts of space that are necessary for capturing the important dynamics of the system, computational time can be saved without losing a significant amount of accuracy.

In addition to the observed spatial heterogeneity in systems where molecules move only by pure diffusion in free space, there can be spatial effects in systems due to interactions with internal lower dimensional structures in the cells. This could for instance be binding to and reactions on curved surfaces or one-dimensional structures such as DNA or microtubules. In Ref. \onlinecite{HeHeLo} a method for simulating interactions with curved surfaces at the microscopic level is developed, and in Ref. \onlinecite{2013arXiv1302.0793M} they develop an FPKMC method for simulating active transport due to chemical potentials on embedded lines of fixed length and with fixed configuration. A mesoscopic method that simulates interactions between molecules in free space and lines is considered in Ref. \onlinecite{it:2012-034} and a mesoscopic method incorporating active transport is suggested in Ref. \onlinecite{HellanderLotstedt2010}.

In this paper we have developed an algorithm that combines single molecule simulations in three-dimensional space with single molecule simulations on embedded curves. Molecules in some general three-dimensional volume, for instance a sphere or a cylinder, are simulated with the GFRD algorithm. They can react with the embedded curves according to the Smoluchowski model. Molecules attached to a curve can also diffuse and react with each other according to the one-dimensional Smoluchowski equation. In addition to pure diffusion, we simulate active transport along the curves by letting the molecules bound to the curves move deterministically in one direction. The curves can have an arbitrary shape. This flexibility is obtained by approximating a continuous curve with non-constant curvature by straight linear segments, and this approximation can be made arbitrarily accurate by making the discretization of the curve finer and finer. It is also known that many one-dimensional structures in cells are highly dynamic \cite{MiKi84}. A microtubule may for instance both move in space as well as grow and shrink. We take those dynamics into account by letting the embedded curves move according to some general random process as well as letting them grow and shrink over time.

The following sections are organized as follows. In Sec. \ref{sim3D1D} we review the GFRD method in one dimension as well as three dimensions. In Sec. \ref{emb_1D} we describe the algorithm for coupling the one-dimensional simulations with the three-dimensional simulations, as well as simulation of moving lines embedded in space. In Sec. \ref{results} we demonstrate the flexibility of the algorithm with four examples. In the first example we verify the accuracy of the method by comparing with mesoscopic simulations and theoretical predictions. In the second example we consider two spirals embedded in a sphere. In the third example we consider several straight lines embedded in a spherical domain. Molecules can bind to the lines and are then moved by active transport along the lines. The lines move in space according to a random process, and we show how the active transport affects the distribution of the molecules in space. In the final example we consider lines that grow and shrink. We show how the number of molecules attached to a line will vary with the total length of the lines. 

All parameters are given in SI-units.

\section{Simulations in three and one dimensions}
\label{sim3D1D}
We consider a system of $N$ molecules, where the molecules can undergo bimolecular or unimolecular reactions. The evolution of the system is modeled by the Smoluchowski equation with a mixed boundary condition at the reaction radius of the molecules. In principle one could write down the Smoluchowski equation for the full system, but due to the high dimension it would not be tractable to solve. In GFRD \cite{ZoWo5a} the main idea is to divide the full $N$-body problem into smaller subproblems consisting of the simulation of single molecules and pairs of molecules. The subdivision is done in two steps. First we find the nearest neighbor of every molecule. Molecules that have each other as nearest neighbors form pairs and the remaining molecules will be updated as single molecules. Secondly we choose a time step $\Delta t$ such that no molecule in a pair is likely to react with any other molecule except its nearest neighbor and molecules that are not in a pair should be unlikely to react with any other molecule. The full $N$-body problem can now be approximated by smaller, independent, subproblems during the time step $\Delta t$.

\subsection{Updating positions of single and pairs of molecules}
\label{sim3D}
The positions of single molecules are updated by sampling from a normal distribution with standard deviation $\sqrt{2D\Delta t}$ in each direction, where $D$ is the diffusion constant of the molecule. A pair of molecules, where the molecules have the positions $\fatx_{1n}$ and $\fatx_{2n}$ at time $t_n$, are updated according to the PDF $p(\fatx_1,\fatx_2,t|\fatx_{1n},\fatx_{2n},t_n)$ obtained by solving the Smoluchowski equation for a pair of molecules, given by:
\begin{align}
\label{eq:smol_full_3D}
p_t = D_1\Delta_{\fatx_1}p+D_2\Delta_{\fatx_2}p,
\end{align}
where $D_i$ and $\fatx_i$ are the diffusion coefficient and position of particle $i$. By making a change of variables as suggested in Ref. \onlinecite{ZoWo5a} this equation can be split into two equations: one for the relative position of the molecules and one for a weighted mean of the positions
\begin{align}
\label{eq:coordinates}
\fatY = \sqrt{D_2/D_1}\fatx_1+\sqrt{D_1/D_2}\fatx_2, \,\, \faty = \fatx_2-\fatx_1.
\end{align}
Now $p$ can be written as $p(\fatx_1,\fatx_2,t|\fatx_{1n},\fatx_{2n},t_n) = p_{\fatY}(\fatY,t|\fatY_n,t_n)p_{\faty}(\faty,t|\faty_n,t_n)$, and the equations become
\begin{align*}
\partial_t p_{\fatY} = D\Delta_{\fatY}p_{\fatY}, \,\, \partial_t p_{\faty} = D\Delta_{\faty}p_{\faty},
\end{align*}
where $D = D_1+D_2$ is the sum of the diffusion coefficients of the molecules. The motion in the $\fatY$ direction is in free space and will be described by a normal distribution, but the equation for the $\faty$ coordinate will have an additional boundary condition at the reaction radius of the molecules to model associations between molecules. By considering the equation in a spherical coordinate system, $\fatr=(r,\theta,\phi)$, the boundary condition becomes:

\begin{align}
\label{eq:3d-bcond}
 4\pi\sigma^2 D\frac{\partial p_{\fatr}}{\partial r}\Big|_{r=\sigma} = k_r p_{\fatr}(r=\sigma, t|r_n,t_n),
\end{align}

where $\sigma$ is the reaction radius of the molecules, $r=\|\fatr\|$ and $k_r$ is the intrinsic association rate. The initial condition is given by 
\begin{align*}
p_{\fatr}(\fatr,t|\fatr_n,t_n)=\delta(\fatr-\fatr_n)
\end{align*}
and $p_{\fatr}(\fatr,t|\fatr_n,t_n)$ vanishes as $r\to\infty$.

This equation can be solved analytically \cite{CarJae,ZoWo5a}, but due to the complexity of the analytical solution we use the method described in Ref. \onlinecite{SHeLo11}. There the full equation is solved in two steps using a first order splitting scheme (or if higher accuracy is needed with a second order splitting scheme). This simplifies the sampling of new positions for the molecules, at the price of an additional error. This error is shown to be small in Ref. \onlinecite{SHeLo11}. The first step of the operator splitting scheme is to solve for the radial part of the equation (in a spherical coordinate system rotated such that $\fatr_n=(r_n,0,0)$), given by:

\begin{align}
\label{eq:3d-radial}
\partial_t p_{r}=D\left(\frac{\partial^2 p_{r}}{\partial^2 r}+\frac{2}{r}\frac{\partial p_{r}}{\partial r}\right),
\end{align}

with the initial condition $p_r(r,t_n|r_n,t_n)=\delta(r-r_n)$ and with the boundary condition given by Eq. \eqref{eq:3d-bcond}. The analytical solution is derived in Ref. \onlinecite{KimShin}. The second step is to solve for the angular part:

\begin{align}
\label{3d-angular}
\partial p_{\theta} = \frac{D}{r^2}\left(\frac{1}{\sin \theta}\frac{\partial}{\partial\theta}\left(\sin\theta\frac{\partial p_{\theta}}{\partial^2\theta}\right)+\frac{1}{\sin^2\theta}\frac{\partial^2 p_{\theta}}{\partial^2\phi}\right),
\end{align}

with initial condition given by $p_{\theta}(\theta,t_n|\theta_n,\phi_n,t_n)=\delta(\theta)/r^2\sin(\theta)$. The solution to Eq. \eqref{3d-angular} will be independent of $\phi$, and the distribution of $\phi$ will therefore be uniform on the interval $[0,2\pi)$. It is possible to solve Eq. \eqref{3d-angular} analytically and the solution is given in the form of an infinite series, where the terms are products of Legendre functions and exponentials. Efficient strategies to compute the analytical solutions of Eqs. \eqref{eq:3d-radial} and \eqref{3d-angular} are discussed in Ref. \onlinecite{HaHeLo}. An efficient method to generate random numbers from the PDF $p_{\theta}$ is suggested in Ref. \onlinecite{CaEkEl10}.

Note that the probability for two molecules to survive until time $t_n+\Delta t$ is given by
\begin{align*}
S(t_n+\Delta t|\fatr_n,t_n) = 1-4\pi\sigma^2 D\int_{t_n}^{t_n+\Delta t}\frac{\partial p_{\fatr}}{\partial r}\Big|_{r=\sigma}d\tau.
\end{align*}
From $S$ we can sample the exact time when the pair of molecules react (cf. Ref. \onlinecite{ZoWo5a}). Thus, the first step is to check whether the two molecules react before $t_n+\Delta t$. If they do, we execute the reaction at the sampled time, and if they do not react we sample new positions for each molecule by

\begin{enumerate}
\item Sample $r_{n+1}$ from the PDF $p_r$.
\item Sample $\theta_{n+1}$ from the PDF $p_{\theta}$ with $r=r_{n+1}$.
\item Sample $\phi_{n+1}$ from a uniform distribution on the interval $[0,2\pi)$.
\item Transform $\fatr_{n+1}$ back to $\faty_{n+1}$. Sample $\fatY_{n+1}$ from a normal distribution and then solve Eq. \eqref{eq:coordinates} to obtain $\fatx_{1(n+1)}$ and $\fatx_{2(n+1)}$.
\end{enumerate}

Also molecules in one dimension can react. For pairs of molecules we solve the one-dimensional Smoluchowski equation
\begin{align}
\displaystyle{\partial_t p_{s}=D\frac{\partial^2 p_{s}}{\partial s^2},}
\label{eq:preq1D}
\end{align}
with the boundary condition at the reaction radius
\begin{align}
D\frac{\partial p_{s}}{\partial s}|_{s=\sigma} = k_r p_{s}(s=\sigma,t|s_n,t_n).
  \label{eq:prcond1D}
\end{align}
to obtain the PDF $p(s,t|s_n,t_n)$ for the distance $s$ at time $t$ between the two molecules, given that the distance was $s_n$ at time $t_n$. The analytical solution to Eq. \eqref{eq:preq1D} with boundary condition given by Eq. \eqref{eq:prcond1D} is derived in Ref. \onlinecite{CarJae}. 

Both molecules in one dimension and molecules in three dimensions can dissociate. Consider a reaction $A\xrightarrow{k_d}B+C$. The time until an $A$-molecule dissociates is sampled from an exponential distribution with mean $k_d$. After a molecule $A$ has dissociated the $B$- and $C$-molecules are placed at a distance equal to the reaction radius.

\subsection{Molecules close to a curved boundary}
\label{sim2D}
Most biologically realistic geometries are fairly complex. Cells can be close to spherical or have a shape reminding of a spheroid. Eukaryotes typically have complex internal two-dimensional structures. One flexible and efficient way to deal with such non-regular surfaces is to discretize the surface with an unstructured mesh.

In this paper we do not consider volumes with internal two-dimensional structures, and are consequently only concerned with the outer boundaries of the volume. The outer boundary is approximated by a triangular surface mesh. The vertex of each triangle is also a center vertex in the dual mesh generated from the primal mesh. Consider a center vertex $\psi_i$. This center vertex is the vertex of several triangles $T_1,\ldots,T_M$ in the mesh, and to each dual element we associate one linear plane defined by the center vertex $\psi_i$ and with normal given by the average of all the normals $\fatn_1,\ldots,\fatn_M$ to the triangles $T_1,\ldots,T_M$. We have now assigned one linear plane to each dual element in the mesh, and these planes are an approximation of the surface. 

Now, assume that a molecule is close to one of the planes approximating the surface. Choose a time step $\Delta t$ small enough such that the interaction between the molecule and the plane can be approximated by the interaction between an infinite flat plane and a molecule during $\Delta t$. Then consider the change of coordinate system where the molecule diffuses parallel to the plane in two directions and along the normal of the plane in the third direction. The molecule can then be updated according to normal diffusion in the first two directions and in the third direction by sampling from a normal distribution and then, in case the molecule ends up on the wrong side of the plane, simply reflecting that coordinate in the plane.

This approach gives great flexibility to the algorithm and enables us to simulate molecules in different complicated geometries. The details of this approach is described in Ref. \onlinecite{HeHeLo}. Another method to simulate molecules close to an absorbing or reflecting boundary has been suggested in Ref. \onlinecite{PhysRevE.66.056701}.

Note that by making a locally planar approximation we can also include reactions between molecules and surfaces by solving the one-dimensional Smoluchowski equation for molecules that are close to the surface, as done in Ref. \onlinecite{HeHeLo}. In this paper, however, we only consider reflecting outer boundaries, but the methodology for simulating interactions with reactive surfaces (outer boundaries as well as internal surfaces) developed in Ref. \onlinecite{HeHeLo} can be combined with the method described in this paper. Furthermore, we could also combine this method with the hybrid method coupling microscopic and mesoscopic simulation algorithms, also developed in Ref. \onlinecite{HeHeLo}.

\section{Embedded one-dimensional manifolds}
\label{emb_1D}
We consider a continuous and differentiable curve $\Gamma $ embedded in space, parametrized by
\begin{align*}
\Gamma(s) = \left(a(s),b(s),c(s)\right),\, 0\leq s\leq 1.
\end{align*}
$\Gamma $ can be approximated by a piecewise linear curve $\gamma$. Denote the linear segments of $\gamma$ by $\gamma_1,...,\gamma_N$. The segment $\gamma_i$ is defined by the points $\fatp_i^1$ and $\fatp_i^2$. The points are chosen such that $\fatp_i^1=\fatp_{i-1}^2$ and thus all segments are connected. The curve $\Gamma$ can be approximated to arbitrary accuracy by letting $N\to\infty$ in such a way that $\|\fatp_i^1-\fatp_i^2\|\to 0$. 

We want to simulate reactions between molecules in space and the curve as well as diffusion and reactions between molecules on the curve. Furthermore, we want to simulate moving curves that grow and shrink.

\subsection{Interactions with the curve}
\label{simWithCurve}
If a molecule $A$ is close to $\gamma$, we will choose a time step $\Delta t$ such that the molecule is unlikely to interact with more than one of the linear segments during that time step, and such that the linear segment can be approximated, to high accuracy, by an infinitely long straight line during $\Delta t$. The problem of simulating a molecule close to an infinitely long straight line can be reduced to a two-dimensional problem, which can be seen by considering the position of the molecule at time $t_n$, $\fatx_n$, in cylindrical coordinates $\fatr=(r,\theta,z)$, with the $z$-coordinate in the direction of the line. The diffusion of the molecule in the $z$-coordinate is now independent of the diffusion in the $r$ and $\theta$ directions and the interaction with the line.

Assume that the closest linear segment to the molecule is $\gamma_i$. Now consider the sphere $S$ with center at $\fatx_n$ and radius $R=\min\{\|\fatx_n-\fatp_i^1\|,\|\fatx_n-\fatp_i^2\|\}$. By choosing the time step $\Delta t$ such that the molecule $A$ is unlikely to diffuse outside of the sphere $S$ during the time step, the requirement that the molecule is unlikely to react with more than one linear segment is satisfied. If the diffusion coefficient of $A$ is $D_A$, the molecule will, on average, diffuse the distance $\sqrt{6D_A\Delta t}$ during a time step $\Delta t$. Thus, we choose the time step to be
\begin{align*}
\Delta t = \frac{R^2}{K\cdot 6D_A},
\end{align*}
where $K$ is a constant chosen sufficiently large in order to ensure that the probability for the molecule to diffuse outside of the sphere is small. Note that if the molecule is really close to $\fatp_i^1$ or $\fatp_i^2$ the time step could become very small, making the simulation too slow. To avoid this problem we have to choose a smallest value that $\Delta t$ can take, thus introducing an error. However, with a carefully chosen smallest value of $\Delta t$ that error will be small. 

The interaction between the molecule and the linear segment can now be approximated by the two-dimensional Smoluchowski equation during $\Delta t$. The two-dimensional Smoluchowski equation for the relative position (in cylindrical coordinates) is given by:
\begin{align}
\label{eq:smol2D_full}
p_t = D\Delta p,
\end{align} 
with the boundary condition
\begin{align}
\label{eq:smol2D_bc}
 2\pi\sigma D\frac{\partial p_{\fatr}}{\partial r}\Big|_{r=\sigma} = k_r p_{\fatr}(r=\sigma, t|\fatr_n,t_n).
\end{align}

Again we solve the equation in two steps with a first order operator splitting scheme. The radial part is given by:
\begin{align}
\label{eq:smol2D}
\partial_t p_{r}=\left(\frac{\partial^2 p_r}{\partial r^2}+\frac{1}{r}\frac{\partial p_r}{\partial r}\right)
\end{align}
with the boundary condition given by Eq. \eqref{eq:smol2D_bc}, the initial condition given by $p_{r}(r,t_n|r_n,t_n)=\delta(r-r_n)$ and with $p_{r}$ vanishing at infinity. The solution of Eq. \eqref{eq:smol2D} with the boundary condition given by Eq. \eqref{eq:smol2D_bc} is solved analytically in Ref. \onlinecite{CarJae}. The analytical solution is
\begin{align*}
p_{r}(r,t|r_0,t_n) = \frac{1}{2\pi}\int_0^{\infty}\exp(-Du^2(t-t_n))C(u,r,k,k_r)C(u,r_0,k,k_r)udu,
\end{align*}
where $k=2\pi\sigma D$ and the function $C$ is defined by
\begin{align*}
C(u,r,k,h) = \frac{J_0(ur)(kuY_1(\sigma u)+hY_0(\sigma u))-Y_0(ur)(kuJ_1(\sigma u)+hJ_0(\sigma u))}{((kuY_1(\sigma u)+hY_0(\sigma u))^2+(kuJ_1(\sigma u)+hJ_0(\sigma u))^2)^{\frac{1}{2}}}.
\end{align*}
Here $J_0$ and $J_1$ are Bessel functions of the first kind and $Y_0$ and $Y_1$ are Bessel functions of the second kind. Despite having the analytical solution available we choose to compute the solution using a finite difference method, as the analytical solution is fairly complex and expensive to evaluate. The solution could also be tabulated as described in Ref. \onlinecite{HaHeLo}. Tabulating the solution is more efficient but has the drawback of introducing an additional error as well as consuming more memory. 

The second step is to solve for $\theta$ and $z$
\begin{align}
\label{eq:smol2D-theta}
\partial p_{\theta z} = D\left(\frac{1}{r^2}\frac{\partial^2 p_{\theta z}}{\partial \theta^2}+\frac{\partial^2 p_{\theta z}}{\partial z^2}\right),
\end{align}
with initial condition
\begin{align*}
p_{\theta z}(\theta,z,t_n)=\frac{\delta(\theta)\delta(z-z_n)}{r}.
\end{align*} 
Again, $p$ vanishes at infinity. The solution of Eq. \eqref{eq:smol2D-theta} is a Gaussian in both the $z$-direction and the $\theta$-direction.

Molecules can also dissociate from the lines. The time until a dissociation is assumed to be exponentially distributed with mean equal to the intrinsic dissociation rate, and after a dissociation from a line the molecule is placed at a distance equal to the reaction radius between the molecule and the line. As we are using a finite difference scheme to compute the solution of Eq. \eqref{eq:smol2D} we are in fact forced to place the dissociating molecule at a slightly larger distance from the line than the reaction radius, in order to have discretization points between the molecule and the line. This introduces a small error, but it can be made arbitrarily small by placing the molecule closer and closer to the line while increasing the number of discretization points used in the finite difference scheme.

\subsection{Diffusion and reactions on moving curves}
\label{simOnCurve}
A new position at time $t+\Delta t$, $\fatx_{n+1}$, of a single molecule moving by pure diffusion, bound to $\Gamma$ at position $\fatx_n$ at time $t$, where $\fatx_n\in\gamma_i$, is computed by sampling a distance $d$ from a one-dimensional normal distribution with standard deviation $\sqrt{2D\Delta t}$ and then updating
\begin{align}
\label{eq:single_on_line}
\fatx_{n+1}=\fatx_n+d\frac{(\fatp_i^2-\fatp_i^1)}{\|\fatp_i^2-\fatp^1_i\|}.
\end{align}
It is possible that $\fatx_{n+1}$ is no longer on the segment $\gamma_i$. If this is the case we project $\fatx_{n+1}$ to the closest point on $\gamma$. The time step $\Delta t$ should thus be chosen such that $\sqrt{2D\Delta t}$ is small in comparison to $\|\fatp_i^2-\fatp_i^2\|$, so that the molecules do not diffuse a long distance compared to the distance of the segments. Since we are projecting the molecule onto the line we are introducing a small error and the effective diffusion will be slightly slower than the exact diffusion process. However, if the angles between the linear segments are small, this error will also be small, and it is therefore important to ensure that the discretization of $\Gamma$ is fine enough to make the angles small enough. Note that the process defined by Eq. \eqref{eq:single_on_line} can be modified to account for other types of motion on curves, such as active transport. For instance, $d$ could be sampled from another distribution than the normal distribution, or be completely deterministic. Molecules on the same line segment can react according to the PDF obtained by solving Eq. \eqref{eq:preq1D} with boundary condition given by Eq. \eqref{eq:prcond1D}.

Many one-dimensional structures in cells are not stationary, but rather move around in space. They can also grow and shrink. We simulate this process with an operator splitting scheme. Choose a time step $\Delta t_{\mathrm{split}}$. First we keep the lines constant in space and propagate the molecules according to the algorithm described above for the full time step $\Delta t_{\mathrm{split}}$. The second step is to move the lines according to some transformation $T(\fatp,t):\mathbb{R}^3\times \mathbb{R}_+\to\mathbb{R}^3$ in space, for the full time step $\Delta t_{\mathrm{split}}$. $T$ can both move the line in space but also change the length of the line. We apply $T$ to all the points, such that the linear segments at time $t+\Delta t_{\mathrm{split}}$ are defined by the pair of points
\begin{align*}
(T(\fatp_1^1,\Delta t_{\mathrm{split}}),T(\fatp_1^2,\Delta t_{\mathrm{split}})),\ldots,(T(\fatp_N^1,\Delta t_{\mathrm{split}}),T(\fatp_N^2,\Delta t_{\mathrm{split}})).
\end{align*}

The position of the molecules are now no longer guaranteed to be on the lines. Each molecule is therefore projected down to the closest point on its respective line. For high accuracy we should therefore choose the time step $\Delta t_{\mathrm{split}}$ such that each molecule is projected only a short distance. After each step it is possible that some molecules on the line that were previously separated by some small distance are now overlapping. If this is the case they are moved apart and placed at a distance equal to their reaction radius. Molecules in space that overlap with the lines are also moved such that they are at a distance from the line equal to their reaction radius. We then repeat this process until the final time $t_f$. If higher accuracy is needed it would be possible to consider higher order splitting schemes, such as Strang splitting \cite{Str}.

The algorithm outlined above is summarized in Algorithm 1. Note that this is only a summarized version of the GFRD algorithm. There are some special cases that require extra care, such as pairs of molecules that are also close to a boundary. For a more detailed description of how to handle such cases the reader is referred to Refs. \onlinecite{ZoWo5a,SHeLo11,TaTNWo10}.

\begin{algorithm}
\caption{GFRD with moving one dimensional structures}
\label{alg:GFRDext}
\begin{algorithmic}[1]
\State Initialize the system at the time $t=0$ and choose a final time $t_f$. Choose a time step $\Delta t_{\mathrm{split}}$ for the operator splitting.
\State Divide the set of all molecules into subsets of one or two molecules such that pairs are formed by molecules that are each others nearest neighbor. Set $t_{\mathrm{loc}}=0$.
\WHILE {$t_{\mathrm{loc}}<\Delta t_{\mathrm{split}}$}
\State {Choose a time step $\Delta t\leq\Delta t_{\mathrm{split}}-t_{\mathrm{loc}}$ such that molecules are unlikely to interact with more than one other molecule during that time step. $\Delta t$ should also be chosen such that no molecule is likely to interact with a two-dimensional boundary or a one-dimensional curve and a neighbor during the time step. To fulfill this requirement certain pairs and single molecules may have to be updated in several smaller steps during $\Delta t$.}
\State {Pairs of molecules are updated according to the PDF obtained by solving the three-dimensional Smoluchowski equation \eqref{eq:smol_full_3D}, as described in Sec. \ref{sim3D}.}
\State {Molecules close to a line are updated according to the PDF obtained by solving the two-dimensional Smoluchowski equation \eqref{eq:smol2D}, as described in Sec. \ref{simWithCurve}.}
\State {Molecules on a line are updated according to the PDF obtained by solving the one-dimensional Smoluchowski equation \eqref{eq:preq1D}, as described in Sec. \ref{sim3D}.}
\State {Molecules that dissociate from the cylinder at time $\tau<\Delta t$ are placed at the reaction radius of the molecule and the cylinder. If a molecule in free space dissociates we place the products at a distance equal to the sum of the reaction radii of the two molecules.}
\State {Update $t_{\mathrm{loc}}=t_{\mathrm{loc}}+\Delta t$.}
\ENDWHILE
\State Update the linear segments according to the transformation $T$.
\State Update the time $t = t+\Delta t_{\mathrm{split}}$.
\State Repeat (2)-(12) until the final time $t_f$.
\end{algorithmic}
\end{algorithm}

\section{Numerical results}
\label{results}
We show the practical use of the algorithm by studying both examples with straight lines through a three-dimensional domain as well as lines with non-constant curvature. We demonstrate how to use the algorithm to simulate active transport by letting molecules attached to a line move with a deterministic speed in a predetermined direction. Finally, we let the lines move in space as well as grow and shrink.

\subsection{Polymer with road blocks embedded in a cylinder}

\begin{figure}
\centering
\includegraphics[scale=0.6]{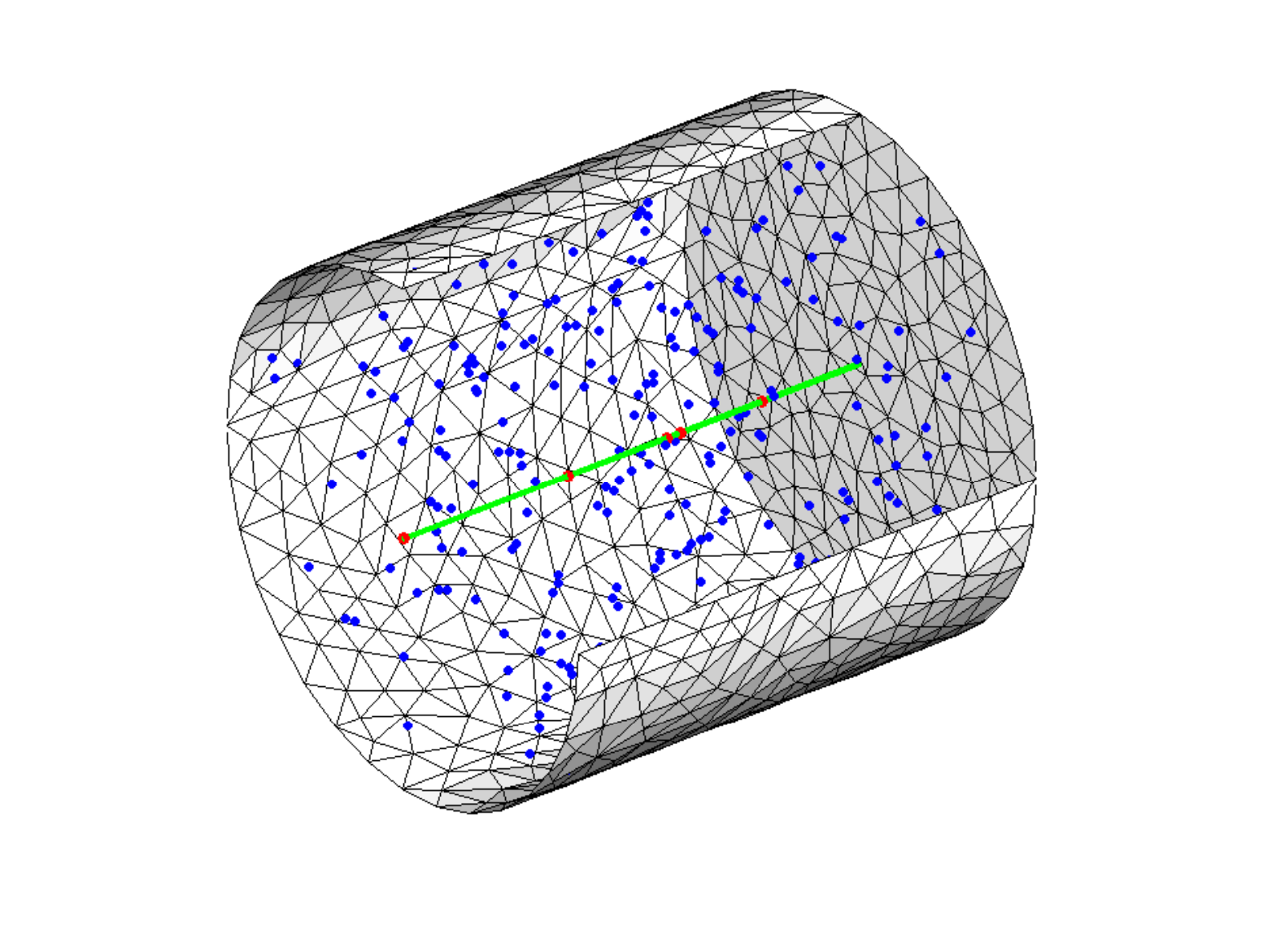}
\caption{The problem setup. Molecules of species $A$ (blue dots) diffuse in the cylindrical volume. They react with the line and transforms into the one-dimensional species $A_{\mathrm{cyl}}$ (red dots). The line (green) extends throughout the domain through the center of the cylinder. The outer boundary is approximated by piecewise linear planes and is reflective.}
\label{fig:linecyl_setup}
\end{figure}

We consider a cylinder with radius $R=10^{-6}$, height $H=2\cdot 10^{-6}$ and aligned along the x-axis. Inside the cylinder we consider a line of length $H$, for instance modeling a polymer, going through the center of the cylinder. The polymer has a reaction radius $r_l=10^{-9}$. Molecules of species $A$ can associate to the polymer to form molecules $A_{\mathrm{cyl}}$. Molecules $A_{\mathrm{cyl}}$ can dissociate from the polymer to become species $A$. Thus, we consider the reaction
\begin{align*}
A\xrightleftharpoons[k_d]{k_r} A_{\mathrm{cyl}}.
\end{align*}
where $A$ is a molecule diffusing in three-dimensional space and $A_{\mathrm{cyl}}$ is a molecule bound to the polymer. The diffusion constant of $A$ is $D_A=10^{-12}$ and the diffusion constant of $A_{\mathrm{cyl}}$ is $D_{A_{\mathrm{cyl}}}=10^{-14}$. The intrinsic reaction rates are chosen to be $k_r=10^{-11}$ and $k_d=50$. The problem setup is shown in Fig. \ref{fig:linecyl_setup}.

\begin{figure}
\centering
\includegraphics[scale=0.60]{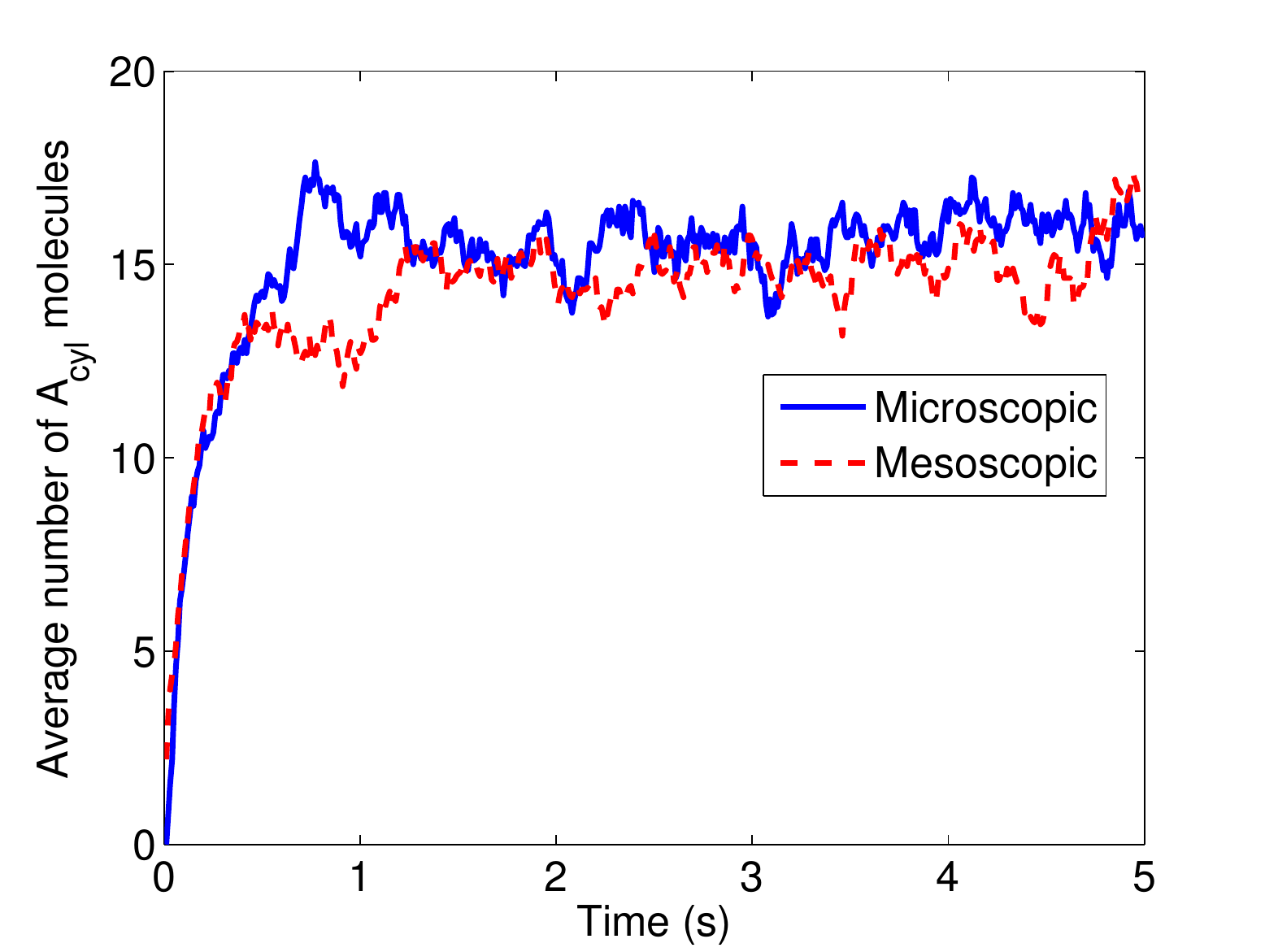}
\caption{The average number of $A_{\mathrm{cyl}}$ at the microscopic level and the mesoscopic level. The mesoscopic simulation agree well with the microscopic simulation. The average is taken over 20 trajectories and we start the simulation with 250 $A$-molecules uniformly distributed in the domain.}
\label{fig:linecyl_res}
\end{figure}

Due to the cylindrical symmetry of the problem, the reactions with the polymer can be reduced to a two-dimensional problem. The average time $T_{\mathrm{bind}}$ until a molecule $A$ reacts with the polymer can thus be computed using a formula for association times on disks, derived in Ref. \onlinecite{FBSE10}. We can therefore compare the microscopic simulation with a corresponding mesoscopic simulation of the reaction
\begin{align*}
A\xrightleftharpoons[k_d^{\mathrm{meso}}]{k_r^{\mathrm{meso}}}A_{\mathrm{cyl}}
\end{align*} 
where $k_r^{\mathrm{meso}}=1/T_{\mathrm{bind}}$ and, due to the detailed balance condition, $k_d^{\mathrm{meso}}=k_d k_r^{\mathrm{meso}}/(A_{\mathrm{disk}} k_r )$. Here $A_{\mathrm{disk}}=\pi R^2$ is the area of the face of the cylinder. In Fig. \ref{fig:linecyl_res} we show that for this simple problem, the microscopic simulations agree well with spatially homogeneous mesoscopic simulations, performed with the SSA. 

Next we consider a more involved example. The geometrical setup is the same as above, but now we let the molecules of species $A_{\mathrm{cyl}}$ react with an operator site at the center of the line, thus considering the line a coarse model of DNA. The operator site is denoted by $B$ and we consider the reactions
\begin{align*}
&A\xrightleftharpoons[k_d]{k_r} A_{\mathrm{cyl}}\\
&A_{\mathrm{cyl}}+B\xrightarrow{} B_{\mathrm{active}}.
\end{align*}
We assume that the operator site is absorbing and of radius $10^{-9}$. On both sides, at the distance $l_{\mathrm{rb}}$ from the operator site, we place a road block. Road blocks are modeled by stationary, non-reactive (but reflective) molecules of radius $10^{-9}$. Thus, molecules binding to the DNA at a distance larger than $l_{\mathrm{rb}}$ from the operator site, will not be able to find it since molecules cannot diffuse through the road blocks. In Ref. \onlinecite{it:2012-034} a similar problem is simulated with a mesoscopic method. The microscopic parameters cannot easily be transformed to corresponding mesoscopic parameters, but the qualitative behaviour can be compared.

\begin{figure}
\centering
\includegraphics[scale=0.6]{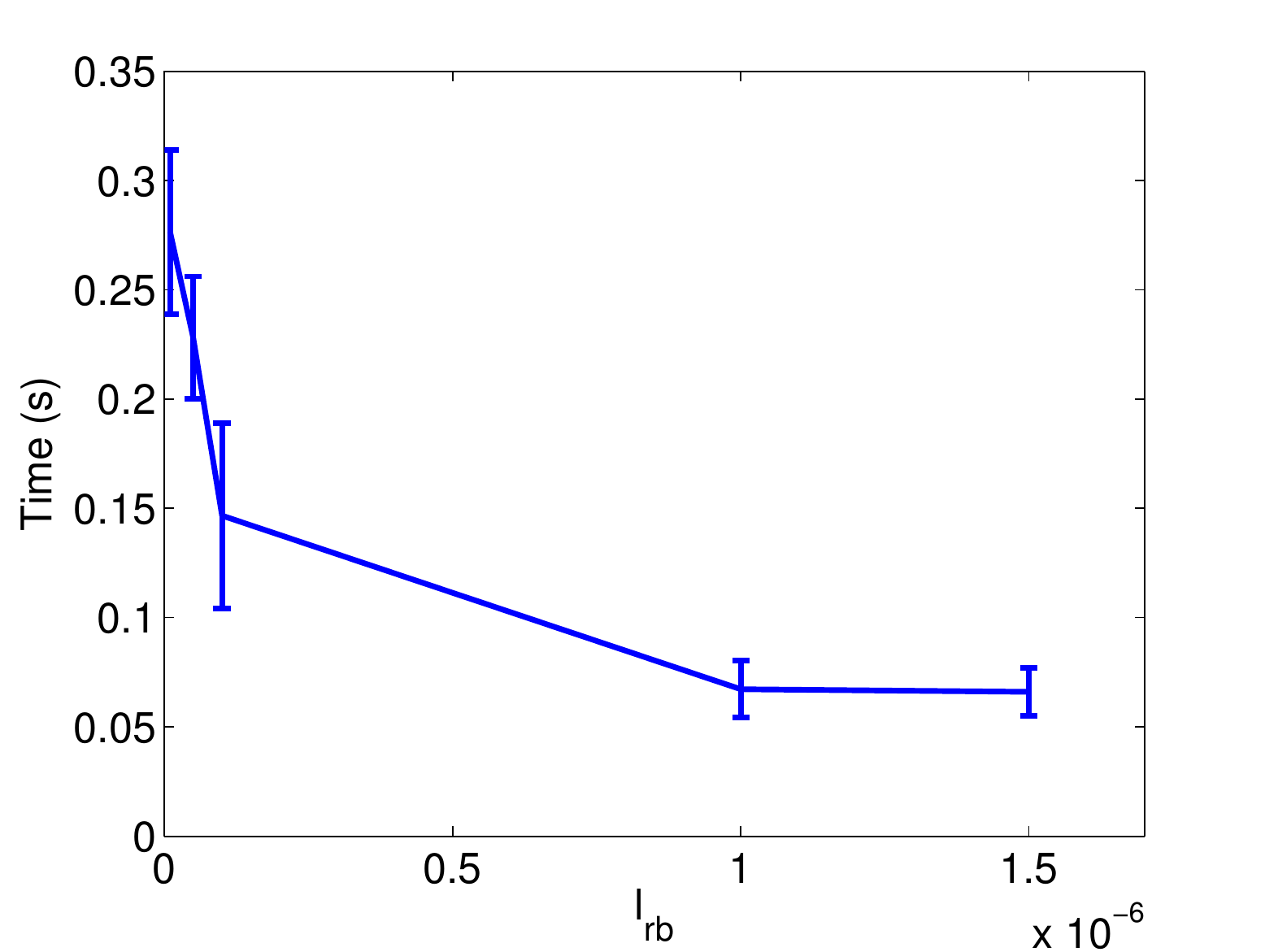}
\caption{The average time until an $A$-molecule binds to the operator site decreases with the distance between the road blocks. As the distance becomes greater than the average sliding distance of an $A_{\mathrm{cyl}}$, the decrease levels off. The parameters in the simulations are $k_r = 10^{-11}$, $k_d = 50$ and $D=10^{-12}$ for all species. The simulation is started with 250 $A$-molecules in free space.}
\label{fig:roadblocks}
\end{figure}

We compute the average time, $T_{\mathrm{active}}$, until an $A_{\mathrm{cyl}}$-molecule binds to the operator site, and this time will depend on the length $l_{\mathrm{rb}}$ between the road blocks. In Fig. \ref{fig:roadblocks} we show that this time decreases with increasing length between the road blocks, as was shown in Refs. \onlinecite{HaLeMaMaBeEl,it:2012-034}. The decrease in the binding time levels off as $l_{\mathrm{rb}}$ becomes greater than the average distance that a molecule slides on the DNA before dissociating. 

\subsection{Spirals embedded in a sphere}

\begin{figure}
\centering
\includegraphics[trim=4cm 2cm 4cm 2cm,clip,scale=0.75]{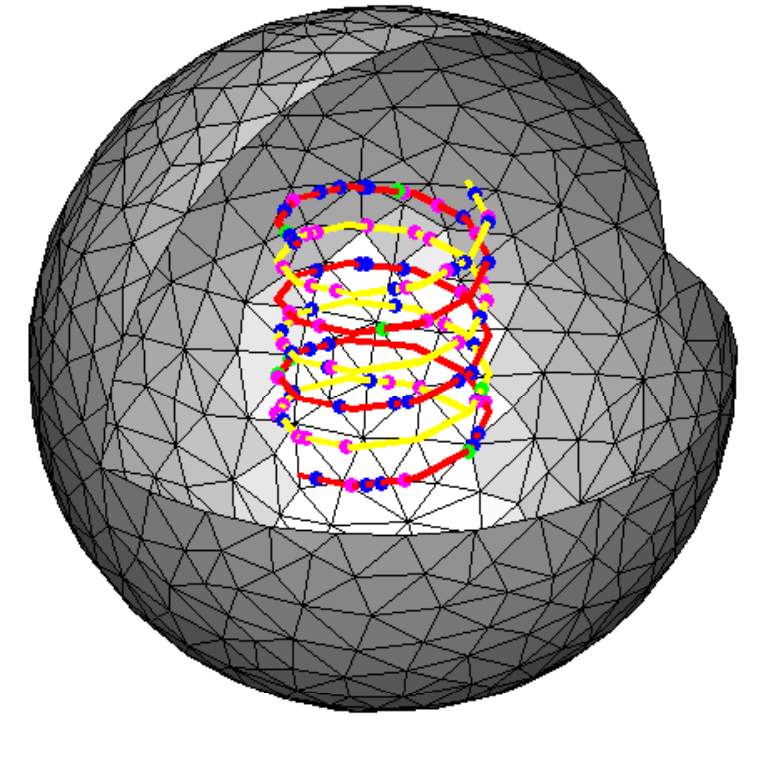}
\caption{We have simulated association to and dissociation from two spirals (yellow and red lines) embedded in a sphere. Molecules in space diffuse and react with the spirals, and molecules on the spirals diffuse and react with each other. Note that the molecules that are free in space are not plotted in the figure above.}
\label{fig:spiral}
\end{figure}

One-dimensional structures in living cells are often not well approximated by straight lines due to their curvature, and they may also have a complicated spatial distribution relative to each other. It is therefore important to be able to simulate molecules reacting with and on lines of non-constant curvature. To demonstrate the flexibility of the method we therefore consider the following system of reactions
\begin{align*}
&A\xrightleftharpoons[k_d^1]{k_r^1} A_{\mathrm{cyl}}, \,\, B\xrightleftharpoons[k_d^2]{k_r^2} B_{\mathrm{cyl}}\\
&A_{\mathrm{cyl}}+B_{\mathrm{cyl}}\xrightleftharpoons[k_d^3]{k_r^3} C_{\mathrm{cyl}},
\end{align*}
where $A_{\mathrm{cyl}}$, $B_{\mathrm{cyl}}$ and $C_{\mathrm{cyl}}$ are molecules bound to one of the lines and $A$ and $B$ are molecules diffusing in three-dimensional space that can react with the lines. We will consider two spirals $\Gamma_1$ and $\Gamma_2$ along the y-axis, parameterized by
\begin{align*}
\Gamma_1(a(s),b(s),c(s))&=(-3\cdot 10^{-7},0,0)+(2.5\cdot 10^{-8}s,r_c\cos(2\pi s),r_c\sin(2\pi s))\\
\Gamma_2(a(s),b(s),c(s))&=(-3\cdot 10^{-7},0,0)+(2.5\cdot 10^{-8}s,-r_c\cos(2\pi s),-r_c\sin(2\pi s)),
\end{align*}
with $0\leq s\leq 3$ and $r_c=3\cdot 10^{-7}$. Each spiral is divided into 30 linear segments. The setup is shown in Fig. \ref{fig:spiral}.

\begin{figure}
\centering
\includegraphics[scale=0.6]{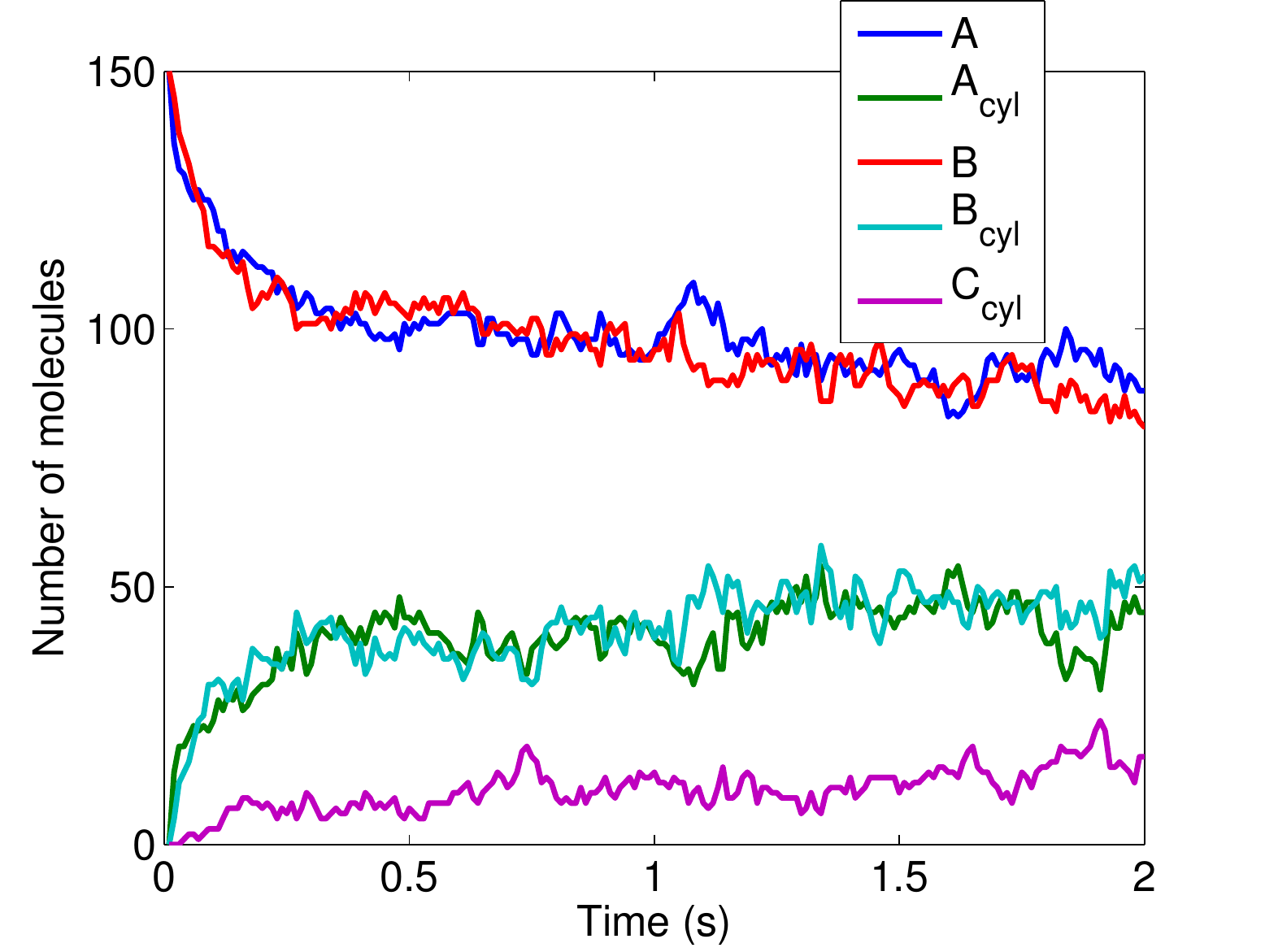}
\caption{The time evolution of the total number of molecules of the different species in the system.}
\label{fig:spiral_res}
\end{figure}

We start with 150 molecules of species $A$ and 150 molecules of species $B$, with initial positions sampled from a uniform distribution. The reactions rates are chosen to be $k_r^1=k_r^2=10^{-11}$, $k_d^1=k_d^2=50$, $k_r^3=10^{-6}$ and $k_d^3=20$. The diffusion rate for the molecules in space is $10^{-12}$ and the diffusion rate for the molecules on the spirals is $10^{-14}$. We simulate the system for $2$ seconds and the evolution of the system is shown in Fig. \ref{fig:spiral_res}.

\subsection{Active transport on moving microtubules}
\label{active_transport}

Molecules in living cells do not only move by normal diffusion but can also be actively transported. For instance, intracellular transport on the cytoskeleton, a structure in the cell composed of fibers such as actin and microtubules, is an important mechanism for directed transport of cargo in the cell \cite{howard}. Cargo attach to the fibers and are then transported by molecular motors: myosins along actin and dynein and kinesins along microtubules. 

We will consider a much simplified model of active transport from the cytosol towards the center of a cell. Molecules bind to microtubules, modeled by straight lines passing through the origin. We assume that the transport is of constant speed in the direction towards the origin and we do not explicitly model the molecular motors. Other more complex and more biologically accurate models of active transport could be incorporated into the method by modifying Eq. \eqref{eq:single_on_line} in a suitable way. The model used in this example is not meant to be a very accurate model of such processes, but is rather meant to serve as an example of the flexibility of the algorithm.

\begin{figure}
\centering
\includegraphics[scale=0.6]{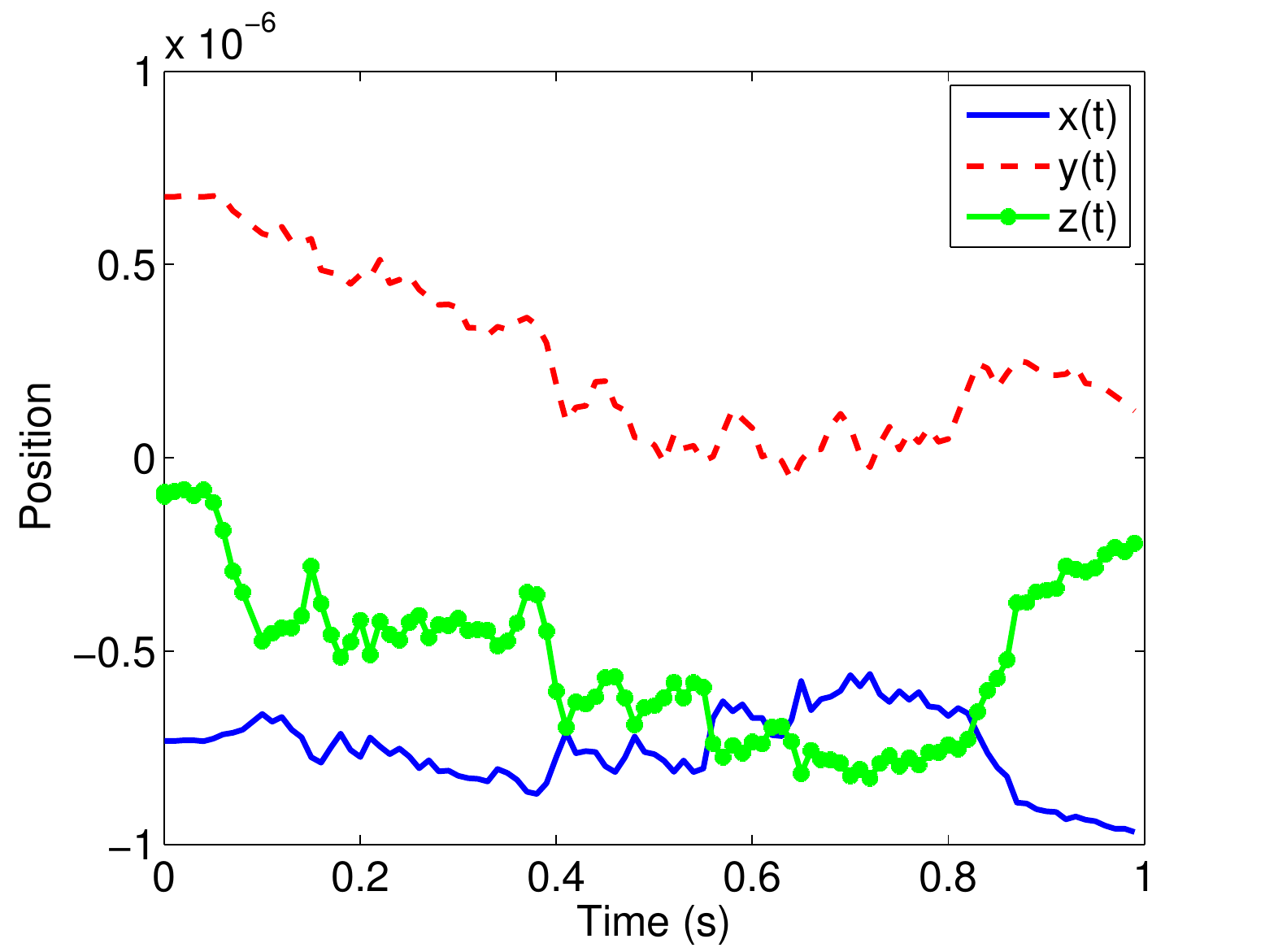}
\caption{The trajectory of one of the microtubules. The position of the $x$-, $y$- and $z$-coordinates of one of the endpoints of the microtubule is plotted. The microtubule is rotated according to Eq. \eqref{eq:line_trans}.}
\label{fig:microtubule_ts}
\end{figure}

\begin{figure}
\centering
\subfigure[$t=0.05s$]{
\includegraphics[trim=4cm 2cm 4cm 2cm,clip,scale=0.65]{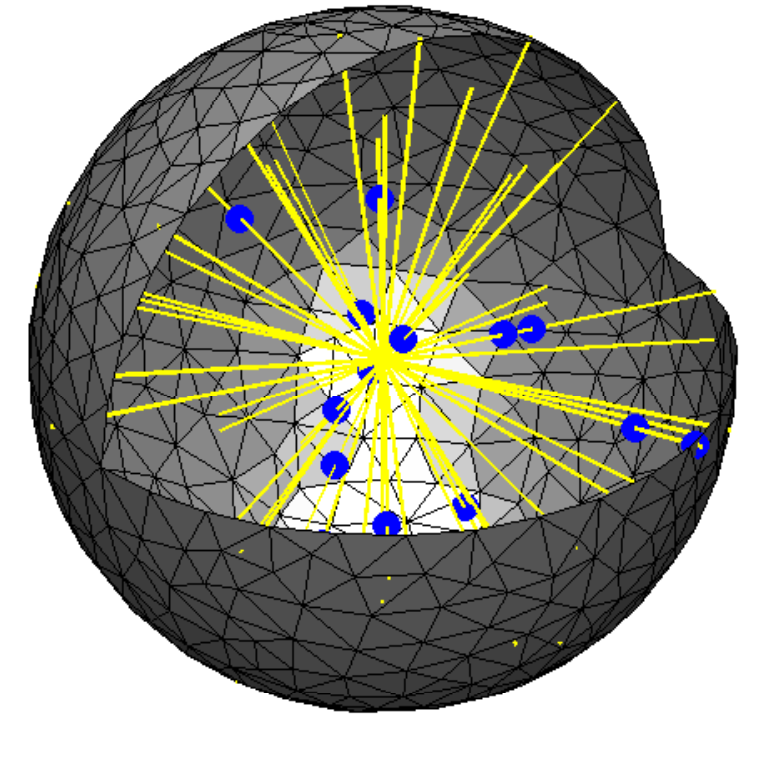}
}
\subfigure[$t=1.0s$]{
\includegraphics[trim=4cm 2cm 4cm 2cm,clip,scale=0.65]{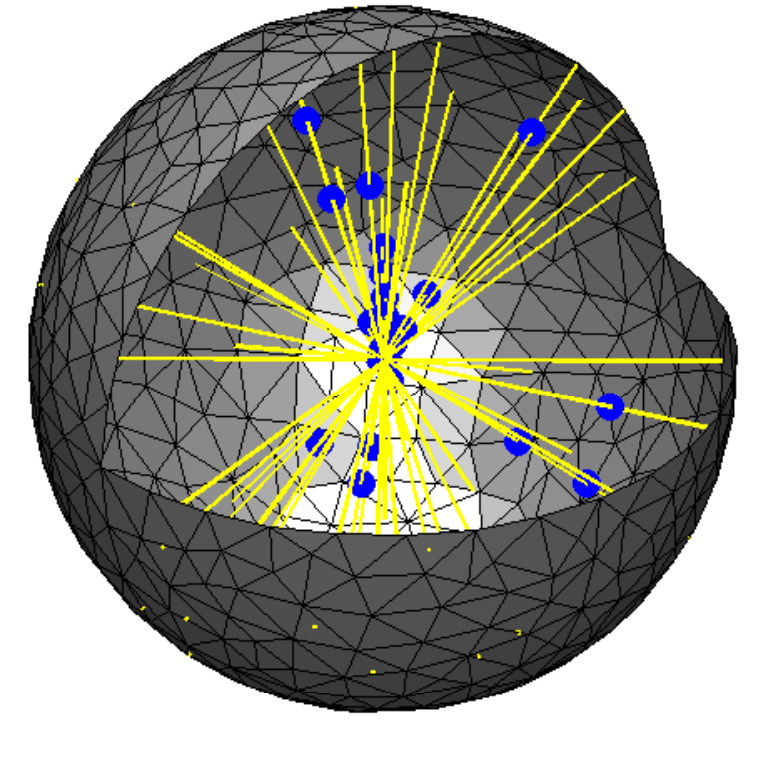}
}
\caption{Lines are yellow and the blue dots are molecules bound to the lines. Molecules in free space are not plotted. Molecules on the lines move by active transport towards the center of the sphere and the lines move in space according to Eq. \eqref{eq:line_trans}.}
\label{fig:active_transport_snapshots}
\end{figure}

We initialize the system by sampling 25 points, $\fatp_1,\ldots,\fatp_{25}$, uniformly on a sphere of radius $R = 10^{-6}$. 25 straight line segments through the origin are then defined by the pairs $(\fatp_i,-\fatp_i)$, $i=1,\ldots,25$. We start with 100 $A$-molecules in space, with initial positions sampled from a uniform distribution. The $A$-molecules can then associate and dissociate to the microtubules through the reaction $A\xrightleftharpoons[k_d]{k_r}A_{\mathrm{cyl}}$. $A$-molecules in space diffuse with the diffusion constant $D=10^{-12}$, but molecules on the microtubules move deterministically towards the center of the sphere with constant speed $\sqrt{2\cdot 10^{-14}\Delta t}/40$. The microtubules are rotated $\Delta\theta$ and $\Delta\phi$ radians around the x- and y-axis per time step $\Delta t$, where
\begin{align}
\label{eq:line_trans}
\begin{split}
\Delta\theta &= 2\pi\Delta t X_1\\
\Delta\phi &= 2\pi\Delta t X_2.
\end{split}
\end{align}
$X_1$ and $X_2$ are random numbers sampled from a normal distribution with standard deviation 1. $\Delta t_{\mathrm{split}}$ is chosen to be 0.01. Thus, the microtubule $l_i$ defined by the points $(\fatp_i,-\fatp_i)$ at time $t$, will be updated at time $t+\Delta t$ by rotating $p_i$ and $-p_i$ around the x- and y-axis $\Delta\theta$ and $\Delta\phi$ degrees, with $\Delta\theta$ and $\Delta\phi$ given by Eq. \eqref{eq:line_trans}. The trajectory of one of the microtubules is plotted in Fig. \ref{fig:microtubule_ts} and two snapshots of a simulation, one at $t_1=0.05s$ and one at the final time $t_f=1.0s$ are shown in Fig. \ref{fig:active_transport_snapshots}.

\begin{figure}
\centering
\includegraphics[scale=0.6]{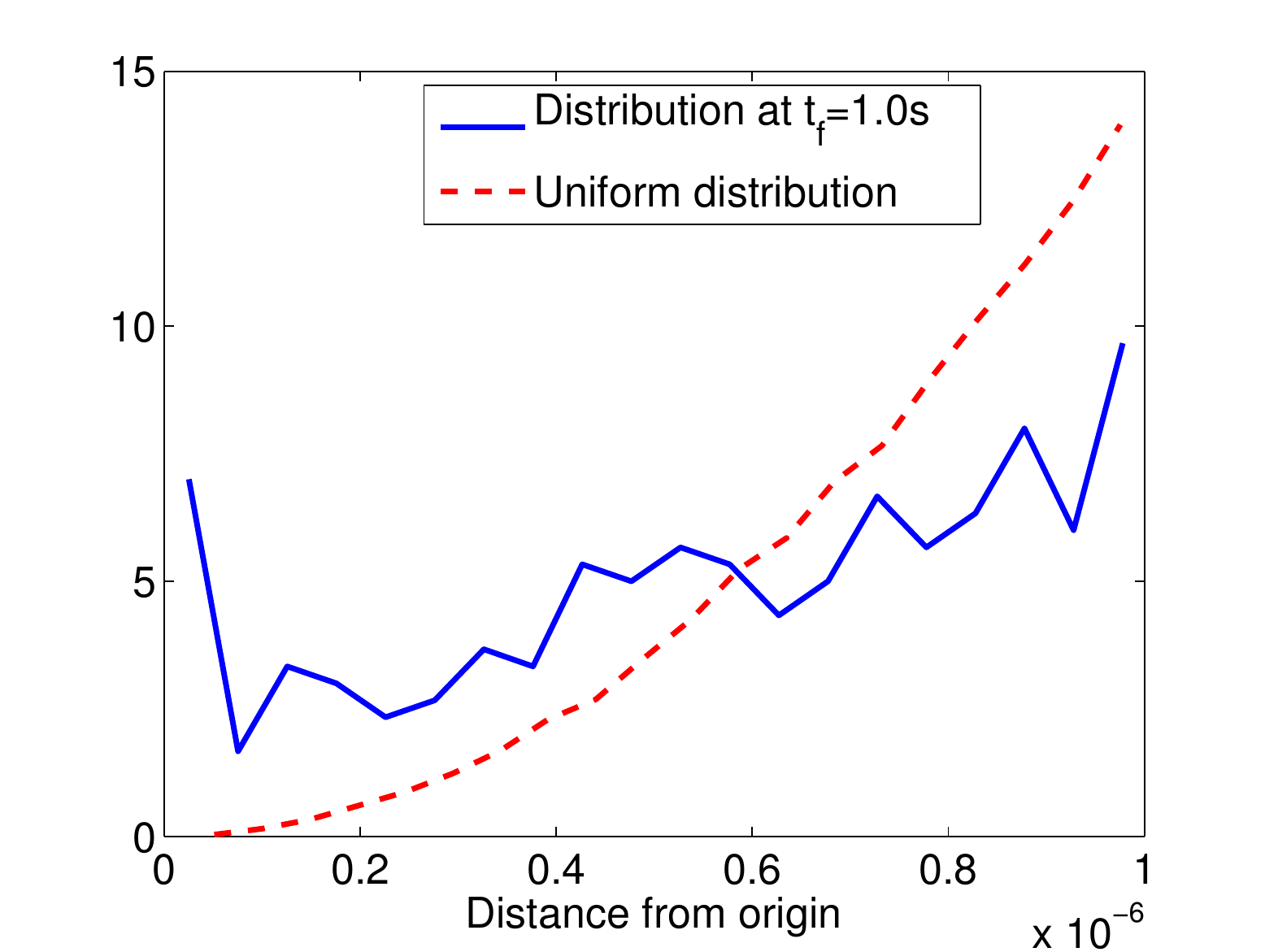}
\caption{The radial distribution of the positions of the molecules. The molecules have a uniform initial distribution in the sphere (dashed line). The distribution at the final time $t_f=1.0s$ is shifted towards the center of the sphere. The distribution at the final time, plotted in the figure above, is based on the mean of three trajectories.}
\label{fig:active_transport_res}
\end{figure}

In Fig. \ref{fig:active_transport_res} the distribution of molecules in the radial direction is plotted. The initial distribution is uniform in space. Due to the transport towards the center, this distribution is shifted and more molecules are close to the center of the sphere at the final time $t_f=1.0s$. The reaction rates in the simulation are $k_r=10^{-11}$ and $k_d=10^3$.


\subsection{Growing and shrinking microtubules}

Microtubules are highly dynamic and can, in addition to moving around in space, both grow and shrink. Therefore, as a final example, we will consider fibers embedded in a sphere and whose lengths can both increase and decrease.

The sphere has a radius $R = 10^{-6}$, with two straight fibers passing through the origin. The fibers are initially of length $l = 2R$. The lines rotate in the same way as in Sec. \ref{active_transport} but will also grow and shrink. After each time step, the length of the line is updated, and molecules are projected down to the closest point on the line. Note that the closest point could in fact be the endpoint of the line if the molecule was close to the end of the line before the length was updated. One possibility is to let such molecules be transformed to molecules in free space, but in this example they will stay bound to the line. Now, assume that the length of the line at time $t_n$ is given by $l(t_n)$. Then the length at $t_n+\Delta t$ is given by
\begin{align*}
l(t_n+\Delta t) = l(t_n)+X,
\end{align*}
where X is a random number sampled from a normal distribution with mean 0 and standard deviation $\sqrt{2D_l\Delta t}$, with the restrictions that $l(t_n+\Delta t)<2R$ and $l(t_n+\Delta t)>10^{-8}$. Thus, if the line is defined by the points $(\fatp_n,-\fatp_n)$ at time $t_n$, it will be defined by the points $l(t_n+\Delta t)/l(t_n)\cdot (\fatp_n,-\fatp_n)$ at time $t_{n+1}$.

Intuitively one would expect that when the length of the fibers decrease, the number of bound molecules will also decrease, under the condition that the fibers grow and shrink at a time scale that is longer than or similar to the time scale for the system to reach steady state. Thus, if we choose parameters such that this condition is fulfilled, we would expect to see a correlation between the number of molecules of species $A_{\mathrm{cyl}}$ and the total length of the two lines.

\begin{figure}
\centering
\includegraphics[scale=0.6]{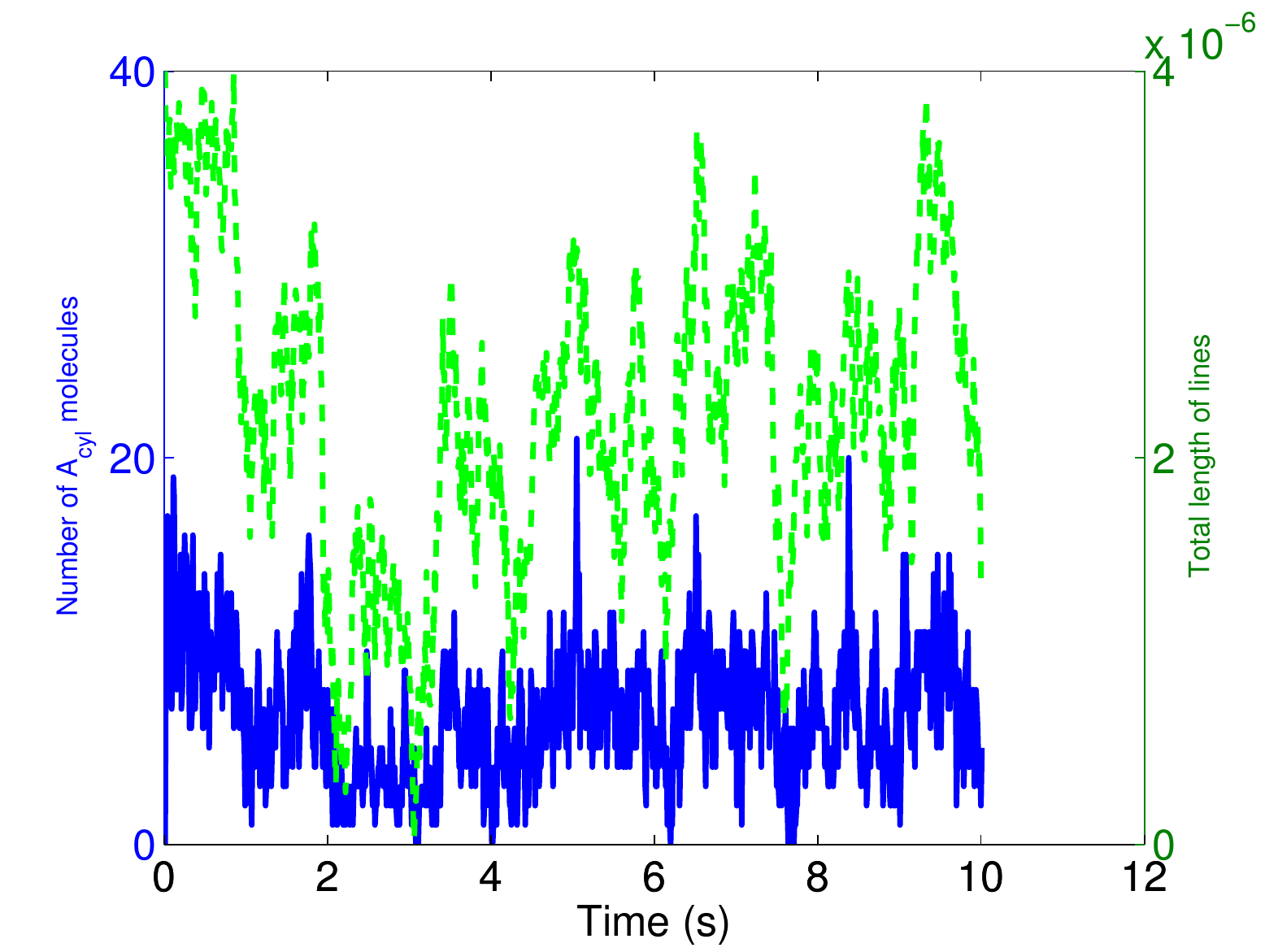}
\caption{The number of $A_{\mathrm{cyl}}$-molecules is plotted together with the sum of the length of the two lines. The two time series are correlated with correlation coefficient 0.713 and a 95\% confidence interval of the correlation coefficient is given by $[0.669,0.752]$.}
\label{fig:growing_sphere_res}
\end{figure}

We start the simulation with $600$ molecules of species $A$ uniformly distributed in space. The diffusion constants are $D_A=10^{-12}$ and $D_{A_{\mathrm{cyl}}} = 10^{-14}$. In Fig. \ref{fig:growing_sphere_res} we plot the total length of the two lines and the total number of molecules bound to a line. In Fig. \ref{fig:growing_sphere_res} a trajectory of the system, with $D_l=10^{-12}$, is plotted together with the sum of the length of the lines. As expected the two trajectories are clearly correlated.

\section{Conclusions}

We have described an algorithm to simulate biochemical systems in complex geometries with embedded and dynamic one-dimensional manifolds. Molecules in space can associate to and dissociate from the embedded curves. Molecules in space and on the lines are simulated using the GFRD algorithm.

The accuracy of the algorithm depends on the discretization of the lines and the time step chosen in the algorithm. The simulation of the molecules and the spatial transformation of the lines are split in time according to a first order operator splitting scheme and the error will therefore depend on the time step chosen for this splitting, as well as the curvature of the lines.

We have demonstrated the flexibility of this method in four numerical examples. First we have shown that the time until a molecule binds to an operator site on DNA depends on the length between two road blocks on each side of the operator site. As the length increases, the time until a molecule binds decreases. The decrease levels off as the distance gets bigger than the average sliding distance of a molecule on the DNA, in agreement with Refs. \onlinecite{HaLeMaMaBeEl,it:2012-034}. In the second example we consider two spirals embedded in a sphere. Molecules in space undergo a reversible reaction with the spirals and molecules on the spirals also react reversibly. In a third example we have considered a simple model of active transport. Molecules in a sphere react reversibly with embedded lines and move deterministically towards the center, while bound. The lines also move in space and we have shown how the radial distribution of the molecules is shifted due to the active transport. In a final example we have simulated lines that are growing and shrinking in space and shown how the dynamics of the lines correlates to the concentration of the molecules in the system.
 
 \section{Acknowledgments}

Andreas Hellander and Per L{\"o}tstedt have read the manuscript carefully and suggested a number of improvements. This work has been supported by the Swedish Research Council and the National Institute of Health under Award no. 1R01EB014877-01. The content is solely the responsibility of the author and does not necessarily represent the official views of the National Institute of Health.

\newcommand{\noopsort}[1]{}

\end{document}